\documentclass[preprint,12pt,3p]{elsarticle}

\usepackage{amssymb}
\usepackage{amsthm}

\usepackage{amsmath}
\usepackage{mathtools}
\usepackage{amsfonts}
\usepackage{mathrsfs}
\usepackage{verbatim}
\usepackage{etoolbox}
\usepackage{hyperref}
\usepackage[none]{hyphenat}

\hypersetup{
    colorlinks=true,
    linkcolor=blue,
    filecolor=magenta,      
    urlcolor=cyan,
}

\newtheorem{thm}{Theorem}[section]

\numberwithin{equation}{section}

\newcommand{\N}{\mathbb{N}}

\newcommand{\R}{\mathbb{R}}

\journal{arXiv}

\begin{document}

\begin{frontmatter}

\title{Wilson's and Wolstenholme's Theorems \tnoteref{label1}}
\tnotetext[label1]{This research did not receive any specific grant from funding agencies in the public, commercial, or not-for-profit sectors.}

\author{Saud Hussein}
\address{Central Washington University, 400 East University Way, Ellensburg, WA 98926}

\ead{saudhussein23@gmail.com}

\begin{abstract}
We provide a proof of Wilson's Theorem and Wolstenholme's Theorem based on a direct approach by Lagrange requiring only basic properties of the primes and the Binomial theorem. The goal is to show how similar the two theorems are by providing the easiest proof possible in a single unified argument.
\end{abstract}


\end{frontmatter}


\section{Introduction}

In 1771, Lagrange \cite{Lagrange} gave the first proof of an interesting property of the prime numbers we now call Wilson's Theorem.

\begin{thm} [Wilson's Theorem] If $p$ is prime, then
\[(p-1)! \equiv -1 \pmod{p}.\]
\end{thm}

Hardy and Wright in their classic book \cite{Hardy} describe Lagrange's approach to this theorem and to an extension called Wolsthenholme's Theorem, first proven by Wolstenholme \cite{Wolstenholme} in 1862.

\begin{thm} [Wolstenholme's Theorem] If $p \geq 5$ is prime, then \[\binom{2p-1}{p-1} \equiv 1 \pmod{p^3}.\]
\end{thm}

In the next section, we lay out Lagrange's simple method in summation notation and in the process see how similar the proofs of Wilson's and Wolstenholme's theorems can be.

\section{Wilson's and Wolstenholme's Theorems}

Let $x \in \R$ and $n \in \N$ odd. Then \begin{align} (x-1)(x-2)\cdots(x-(n-1)) = \sum_{k=0}^{n-1} (-1)^k P_k x^{n-1-k} \label{sym} \end{align} where \[P_k = \sum_{1\leq a_1<a_2<\cdots<a_k\leq n-1}a_1a_2\cdots a_k, \quad a_i\in \N.\]

Note $P_{n-1} = (n-1)!$ and we take $P_0 = 1$, $P_i = 0$ for $i < 0$ and $i\geq n$. Now, multiply both sides of \eqref{sym} by $x$ and make the change of variables $x=y-1$, \begin{align*} x(x-1)(x-2)\cdots(x-(n-1)) &= \sum_{k=0}^n (-1)^k P_k x^{n-k}\\
(y-1)(y-2)\cdots(y-(n-1))(y-n) &= \sum_{k=0}^n (-1)^k P_k (y-1)^{n-k}\\
\left(\sum_{k=0}^{n-1} (-1)^k P_k y^{n-1-k}\right)(y-n) &= \sum_{k=0}^n (-1)^k P_k \left(\sum_{t=0}^{n-k}(-1)^t\binom{n-k}{t} y^{n-k-t}\right)\\
\sum_{k=0}^n (-1)^k (nP_{k-1}+P_k)y^{n-k} &= \sum_{k=0}^n P_k \left(\sum_{m=k}^n(-1)^m\binom{n-k}{m-k} y^{n-m}\right)\\
\sum_{m=0}^n (-1)^m (nP_{m-1}+P_m)y^{n-m} &= \sum_{m=0}^n(-1)^m\left(\sum_{k=0}^m P_k\binom{n-k}{m-k}\right) y^{n-m}.
\end{align*}

In the third equation, we use \eqref{sym} with the arbitrary variable $y$ on the left hand side and the Binomial theorem on the right hand side. In the fourth line, we use a change of variables $m=k+t$ on the right hand side and then switch the summation order to get the final result. Comparing the coefficients on the left and right side, we see that \[nP_{m-1}+P_m = \sum_{k=0}^m P_k\binom{n-k}{m-k},\] so \[(m-1)P_{m-1} = \sum_{k=0}^{m-2} P_k\binom{n-k}{m-k} \quad \text{for} \quad 2\leq m\leq n.\]

Now, assume $n=p$ is prime. Then clearly $p|P_1$ which implies $p|P_2$ which together then imply $p|P_3$ and so on, thus $p|P_i$ for $i=1,2,\dots,p-2$. The case $m=p$ gives us \[(p-1)P_{p-1} = \sum_{k=0}^{p-2} P_k = 1 + \sum_{k=1}^{p-2} P_k\] and since $P_{p-1}=(p-1)!$, we have \[(p-1)!\equiv -1 \pmod{p}.\] The congruence also easily holds for $p=2$ and so we have established Wilson's Theorem.

Since \begin{align} \binom{2p-1}{p-1} = \frac{(p+1)(p+2)\cdots(p+p-1)}{(p-1)!} = \frac{\sum\limits_{k=0}^{p-1} P_k p^{p-1-k}}{(p-1)!} = 1+\frac{\sum\limits_{k=0}^{p-2} P_k p^{p-1-k}}{(p-1)!}, \label{wol} \end{align} then $p|P_i$ for $i=1,2,\dots,p-2$ implies \[\binom{2p-1}{p-1}\equiv 1 \pmod{p^2} \quad \text{for every prime} \quad p\geq 3.\]

Now, let $x=n=p$ be a prime in \eqref{sym}. Then \begin{align*} (p-1)(p-2)\cdots(p-(p-1)) &= \sum_{k=0}^{p-1} (-1)^k P_k p^{p-1-k}\\
(p-1)!&= (p-1)! + \sum_{k=0}^{p-2} (-1)^k P_k p^{p-1-k},
\end{align*}
so \[P_{p-2}= \sum_{k=0}^{p-3} (-1)^k P_k p^{p-2-k}.\] Once again, since $p|P_i$ for $i=1,2,\dots,p-2$, then $p^2|P_{p-2}$ for every prime $p\geq 5$. Therefore by \eqref{wol}, \[\binom{2p-1}{p-1}\equiv 1 \pmod{p^3} \quad \text{for every prime} \quad p\geq 5,\] which is Wolstenholme's Theorem.

We can also easily show the converse of Wilson's Theorem holds. Assume \[(n-1)! \equiv -1 \pmod{n}\] and $n$ is not prime. Since $n$ is not prime, there exists a divisor $d \in \{2,3,4,\dots,n-1\}$ of $n$. But this means $d|(n-1)!$ so $d$ does not divide $(n-1)!+1$, contradicting our assumption. Thus $n$ is prime.


\bibliographystyle{elsarticle-harv}

\bibliography{biblio.bib}

\end{document}